\documentclass[11pt,epsf]{article}
\usepackage[dvips]{color}
\usepackage{indentfirst}
\newcommand{\bed}{\begin{displaymath}}
\newcommand{\eed}{\end{displaymath}}
\newcommand{\bea}{\bed\begin{array}{rl}}
\newcommand{\eea}{\end{array}\eed}
\newcommand{\barray}{\begin{array}{ll}}
\newcommand{\earray}{\end{array}}

\usepackage{color,latexsym,amsfonts,amssymb,amsmath,amsthm}
\usepackage{colordvi,amsfonts,amsmath,epsfig,subfigure,cite}
\usepackage{graphicx}
\usepackage{caption}
\newtheorem{theorem}{Theorem}[section]
\newtheorem{lemma}[theorem]{Lemma}
\newtheorem{remark}{Remark}[section]

\newtheorem{assumption}{Assumption}[section]
\newtheorem{definition}{Definition}[section]

\begin{document}

\title{\Large\bf{Stability of analytical solutions and convergence of numerical methods for non-linear stochastic pantograph differential equations}}
\author{M.H.Song,\thanks{Corresponding author. Email: songmh@lsec.cc.ac.cn,~luyulan2013@163.com,~mzliu@hit.edu.cn
\vspace{6pt}}~~Y.L.Lu,~~M.Z.Liu\\\vspace{6pt}{\em{Department of Mathematics,
Harbin Institute of~ Technology, Harbin, China, 150001}}\thanks{This work is supported by the NSF of P.R. China (No.11071050)}} \maketitle

\begin{abstract}
\indent
 In this paper, we study the polynomial stability of analytical solution and convergence of the semi-implicit Euler method for non-linear stochastic pantograph differential equations. Firstly, the sufficient conditions for solutions to grow at a polynomial rate in the sense of mean-square and almost surely are obtained. Secondly, the consistence and convergence of this method are proved. Furthermore, the orders of consistence (in the sense of average and mean-square) and convergence are given, respectively.

{\bf keywords.} non-linear stochastic pantograph differential equations; polynomial stability; semi-implicit Euler method; consistence; convergence

\end{abstract}

\baselineskip=0.9\normalbaselineskip \vspace{-3pt}
\baselineskip=\normalbaselineskip

\section{Introduction}
Stochastic pantograph differential equations(SPDEs) arise widely in control, biology, neural network, and finance, etc. Asymptotic stability of analytical solution has received  considerable attention in literature for both deterministic and stochastic functional differential equations. Especially, plenty of literature on stability exist with non-exponential rates decay of the solutions(see \cite{cd}\cite{kx}\cite{mz}\cite{ab}). One important non-exponential rates of decay is polynomially asymptotic stability, which means that the rate of decay is controlled by a polynomial function in mean-square or almost surely sense. This type of stability has been studied in \cite{jg1},\cite{jg2} and \cite{jg3}. Buckwar and Appleby consider the polynomial stability of one dimensional linear stochastic pantograph differential equation in \cite{eb}, where the sufficient conditions of polynomially asymptotic property are given.

The convergence of numerical method is another crucial property of stochastic differential equations. Recently, many researchers devoted to the stochastic delay differential equations. Mao Wei\cite{mw} gave the sufficient conditions of convergence with semi-implicit Euler method for variable delay differential equations driven by Poison random jump measure. Fan studied the approximate solution of linear stochastic pantograph differential equations with Razumikhin technique in \cite{fan}. Baker and Buckwar, in \cite{bb}, also investigated the linear stochastic pantograph differential equations, and the sufficient conditions of convergence of semi-implicit Euler method are obtained.

This article investigates non-linear stochastic pantograph differential equations. It is organised as follows. In section 2, we introduce necessary notations and results on pantograph differential equations. The sufficient conditions of polynomial stability of analytical solutions in the sense of mean-square and almost surely are given in section 3. Finally, both the consistence and convergence of the semi-implicit Euler method are proved, and the orders are achieved.
\section{Preliminary notations and properties of pantograph differential equations}
Throughout this paper, unless otherwise specified, the following notations are used. If $x,y$ are vectors, the inner product of $x,y$ is denoted by $\langle x,y \rangle $=$x^Ty$, and $|x|$ denotes the Euclidean norm of $x\in R^d$. We use $A^T$ to denote the transpose of $A$, if $A$ is a vector or matrix. And if $A$ is a matrix, its trace norm is denoted by $|A|=\sqrt{traceA^TA}$.

Let $(\Omega,P,\mathcal{F})$ be a complete probability space with a filtration $\{\mathcal{F}_t\}_{t\geq0}$ satisfying the usual conditions. Let $0<q<1$, $x(t)=(x_1(t),x_2(t),\cdots,x_d(t))^T\in R^d$ and $B(t)=(B_1(t),B_2(t),\cdots,B_m(t))^T$ be a $m$-dimensional Brownian motion. Suppose $h(x)$ be a function, the derivative of $h(x)$ is defined as $D^+h(x)$, i.e

\begin{equation}\label{2.1}
  D^+h(x)=\limsup\limits_{\delta\rightarrow0}\frac{h(x+\delta)-h(x)}{\delta}.
\end{equation}

Before studying the stochastic pantograph differential equations, the properties of deterministic pantograph differential equations are firstly introduced. The equation has the following form:

\begin{equation}\label{2.2}
\left\{
\begin{array}{l}\
x^{\prime}(t)=\bar{a}x(t)+\bar{b}x(qt)
\\x(0)=x_0.
\\
\end{array}
\right.
\end{equation}

\begin{lemma}\cite{eb} Assume $x(t)$ be the solution of (2.2), $x_0>0$, if $\bar{a}<0$, then there exists a constant number $C$ such that

\begin{equation}\label{2.3}
 \limsup\limits_{t\rightarrow\infty}\frac{|x(t)|}{t^\alpha}=C|x_0|,
\end{equation}
where $\alpha\in R$ and $\bar{a}+\bar{b}q^\alpha=0$.
\end{lemma}

\begin{remark} ~If there exists a constant number $C>0$, such that

\begin{equation}\label{2.4}
  |x(t)|\leq C|x(0)|t^\alpha , ~t\geq0,
\end{equation}
then (2.3) holds.
\end{remark}
\begin{lemma}\cite{eb} Assume $x(t)$ be the solution of equation (2.2), $x_0>0$. If $\bar{b}>0$, $p(t):R^{+}\rightarrow R^+$ is a non-negative continuous function satisfying
\begin{equation}\label{2.5}
 D^+p(t)\leq\bar{a}p(t)+\bar{b}p(qt), ~t\geq0,
\end{equation}
and $0<p(0)\leq x(0)$, then $p(t)\leq x(t)$ for all $t\geq0$. Furthermore, if $\bar{a}<0$, and $p(t)$ satisfys (2.5), then there exists a constant $C>0$ such that
\begin{equation}\label{2.6}
  p(t)\leq Cp(0)t^\alpha , t\geq0,
\end{equation}
where $\bar{a}+\bar{b}q^{\alpha}=0$.
\end{lemma}

\section{Polynomial stability of analytical solutions for non-linear stochastic pantograph differential equations}
In this section, we consider the following equation
\begin{equation}\label{3.1}
\left\{
\begin{array}{l}\
   dx(t)=f(t,x(t),x(qt))dt+g(t,x(t),x(qt))dB(t),~0\leq t\leq T
\\
   x(0)=x_0.
\\
\end{array}
\right.
\end{equation}
where  $f:[0,T]\times R^d\times R^d \rightarrow R^d$, $g:[0,T]\times R^d\times R^d\rightarrow R^{d\times m}$, $f(t,0,0)=0$, $g(t,0,0)=0$, and $E|x_0|^2<\infty$. It is easy to see there exists zero solution for (\ref{3.1}).
\begin{definition}
The zero solution for (\ref{3.1}) is said to be mean-square polynomial stable, if there exists a constant number $\alpha<0$, such that
\begin{equation}\label{3.2}
\limsup\limits_{t\rightarrow\infty}\frac{logE|x(t)|^2}{logt}\leq\alpha,
\end{equation}
where $x(t)$ is the solution for (\ref{3.1}) with any initial value $x(0)=x_0$.
\end{definition}
\begin{definition}
The zero solution of (3.1) is said to be almost surely polynomial stable, if there exists a constant number $\alpha<0$, such that
\begin{equation}\label{3.3}
\limsup\limits_{t\rightarrow\infty}\frac{log|x(t)|}{logt}\leq\alpha~~a.s.,
\end{equation}
where $x(t)$ is the solution of (\ref{3.1}) with any initial value $x(0)=x_0$.
\end{definition}
\begin{assumption} ~Assume there exist real numbers $a $, $b>0, c>0 ,d>0$, such that the coefficients $f$ and $g$ satisfy
\begin{equation*}
 \langle x_1-x_2,f(t,x_1,y)-f(t,x_2,y)\rangle\leq a|x_1-x_2|^2,~~\forall x_1,x_2,y\in R^d, t\geq0;
\end{equation*}
\begin{equation*}
  |f(t,x,y_1)-f(t,x,y_2)|\leq b|y_1-y_2|,~~\forall x,y_1,y_2\in R^d, t\geq0;
\end{equation*}
\begin{equation*}
  |g(t,x_1,y)-g(t,x_2,y)|\leq c|x_1-x_2|,~~\forall x_1,x_2,y\in R^d, t\geq0;
\end{equation*}
\begin{equation*}
  |g(t,x,y_1)-g(t,x,y_2)|\leq d|y_1-y_2|,~~\forall x,y_1,y_2\in R^d, t\geq0.
\end{equation*}
\end{assumption}

According to $f(t,0,0)=0$, $g(t,0,0)=0$ and assumption \ref{3.3}, we can estimate
$$\langle x,f(t,x,y)\rangle\leq (a+\frac{b}{2})|x|^2+\frac{b}{2}|y|^2,~~~|g(t,x,y)|^2\leq 2c^2|x|^2+2d^2|y|^2.$$
\begin{theorem} Suppose that $x(t)$ is the solution of (\ref{3.1}), if the coefficients $f,g$ satisfy assumption \ref{3.3}, and $a+b+c^2+d^2<0$, then the zero solution of (\ref{3.1}) is mean-square polynomial stable.
\end{theorem}

{\bf Proof.} According to definition 3.1, we just need to prove that there exist constant number $C_1$ and $\alpha<0$, such that $E|x(t)|^2\leq C_1t^\alpha$.

It\^{o} formula shows that
\begin{equation}
  E(|x(t)|^2)=E(|x_0|^2)+E(\int^t_0[2\langle x(s),f(s,x(s),x(qs))\rangle+|g(s,x(s),x(qs))|^2]{\rm d}s).
\end{equation}
~~~~Let $Y(t)=E(|x(t)|^2)$, then for any $t\geq0$, $t+h\geq0$, we have
\begin{equation*}
\begin{split}
 Y(t+h)-Y(t)&\leq 2E(\int^{t+h}_t[(a+\frac{1}{2}b)|x(s)|^2+\frac{1}{2}b|x(qs)|^2]{\rm d}s)+2E(\int^{t+h}_t[c^2|x(s)|^2+d^2|x(qs)|^2]{\rm d}s)\\
 &\leq (2a+b+2c^2)\int^{t+h}_tY(s){\rm d}s+(b+2d^2)\int^{t+h}_tY(qs){\rm d}s.
 \end{split}
\end{equation*}
Due to
\begin{equation*}
\begin{split}
  Y(t)=\limsup\limits_{h\rightarrow0}\frac{\int^{t+h}_0Y(s){\rm d}s-\int^t_0Y(s){\rm d}s}{h},\\
  Y(qt)=\limsup\limits_{h\rightarrow0}\frac{\int^{t+h}_0Y(qs){\rm d}s-\int^t_0Y(qs){\rm d}s}{h},
 \end{split}
\end{equation*}
thus
\begin{equation}
 D^+Y(t)\leq (2a+b+2c^2)Y(t)+(b+2d^2)Y(qt).
\end{equation}
Note that $2a+b+2c^2<0$,~$b+2d^2>0$, by Lemma 2.3, there exists $C_1$ and $\alpha\in R$, such that
\begin{equation}
  E(|x(t)|^2)=Y(t)\leq C_1Y(0)t^\alpha=C_1E(|x_0|^2)t^\alpha,
\end{equation}
 where $\alpha$ satisfies $2a+b+2c^2+(b+2d^2)q^\alpha=0$. According to $a+b+c^2+d^2<0$, we can know $\alpha<0$, the theorem is proved.

\qed
\begin{theorem} Suppose that $x(t)$ is the solution of (\ref{3.1}), if the coefficients $f,g$ satisfy assumption \ref{3.3}, and $2a+b+2c^2+(b+2d^2)/q<0$, then the zero solution of (\ref{3.1}) is almost surely polynomial stable.
\end{theorem}
{\bf Proof.} According to $2a+b+2c^2+(b+2d^2)/q<0$, we know that $a+b+c^2+d^2<0$ and $E(|x(t)|^2)\leq C_1E(|x(0)|^2)t^\alpha$, where $\alpha<-1$.

By It\^{o} formula and assumption \ref{3.3}, one can show that for any $n-1\leq t\leq n$,
\begin{equation}\label{3.7}
\begin{split}
  E(\sup\limits _{n-1\leq t\leq n}|x(t)|^2)&\leq
  E(|x(n-1)|^2)+(b+2d^2)E\int^n_{n-1}|x(qs)|^2{\rm d}s+2c^2E\int^n_{n-1}|x(s)|^2{\rm d}s
  \\
  &~~~+2E(\sup\limits_{n-1\leq t\leq n}\int^t_{n-1}x^T(s)\cdot g(s,x(s),x(qs)){\rm d}B(s)).
    \end{split}
\end{equation}
According to Burholder-Davis-Gundy inequations, it is easy to show that
\begin{equation}\label{3.11}
\begin{split}
&~~~~2E(\sup\limits_{n-1\leq t\leq n}\int^t_{n-1} x(s)^{T}\cdot g(s,x(s),x(qs)){\rm d}B(s))
\\
&\leq 8\sqrt{2}E(\int^n_{n-1}| x(s)|^2|g(s,x(s),x(qs))|^2{\rm d}s)^\frac{1}{2}
\\
&\leq 8\sqrt{2}E(\sup\limits_{n-1\leq t\leq n}| x(t)|^2\int^n_{n-1}|g(s,x(s),x(qs))|^2{\rm d}s)^\frac{1}{2}
\\
&\leq E(\frac{1}{2}\sup\limits_{n-1\leq t\leq n}| x(t)|^2+64\int^n_{n-1}|g(s,x(s),x(qs))|^2{\rm d}s)
\\
&\leq \frac{1}{2}E(\sup\limits_{n-1\leq t\leq n}| x(t)|^2)+64E(\int^n_{n-1}(c^2|x(s)|^2+d^2|x(qs)|^2){\rm d}s).
\end{split}
\end{equation}

Substituting (\ref{3.11}) into (\ref{3.7}), then
\begin{equation*}
\begin{split}
  E(\sup\limits_{n-1\leq t\leq n}|x(t)|^2)\leq&
  2C_1E(|x_0|^2)(n-1)^\alpha +132c^2C_1E(|x_0|^2)\int^n_{n-1}s^\alpha{\rm d}s
  \\
  &+(2b+132d^2)C_1E(|x_0|^2)q^\alpha\int^n_{n-1}s^\alpha{\rm d}s.
  \end{split}
\end{equation*}
It is easy to know $\int_{n-1}^{n}s^{\alpha}{\rm d}s\leq (n-1)^{\alpha}\max(1,2^{\alpha})$. So
\begin{equation*}
  E(\sup\limits_{n-1\leq t\leq n}|x(t)|^2)\leq \tilde{C}(n-1)^\alpha.
\end{equation*}
where $\widetilde{C}=[2+132c^2+(2b+132d^2)q^\alpha]C_1E(|x_0|^2)\max(1,2^\alpha)$.
By Markov's inequations, for any $\varepsilon>0$, it is not difficult to show
\begin{equation*}
\begin{split}
  &P(\sup\limits_{n-1\leq t\leq n}|x(t)|^2(n-1)^{-1-\alpha-\varepsilon}\geq\gamma)
  \\
\leq& \frac{1}{\gamma}\frac{1}{(n-1)^{1+\varepsilon}}\frac{1}{(n-1)^\alpha}E(\sup\limits_{n-1\leq t\leq n}|x(t)|^2)
  \\
\leq& \frac{1}{\gamma}\frac{1}{(n-1)^{1+\varepsilon}}\widetilde{C}.
\end{split}
\end{equation*}
By using Borel-Cantelli lemma, the following limit can be achieved
\begin{equation*}
  \limsup\limits_{n\rightarrow\infty}\sup\limits_{n-1\leq t\leq n}|x(t)|^2(n-1)^{-1-\alpha-\varepsilon}=0~~~a.s..
\end{equation*}
Note that for any $t>0$, there exists $n(t)$ such that $n(t)-1\leq t\leq n(t)$, and
\begin{equation*}
  \lim\limits_{t\rightarrow\infty}\frac{t}{n(t)-1}=1.
  \end{equation*}
 Hence
  \begin{equation*}
\begin{split}
  \limsup\limits_{t\rightarrow\infty}|x(t)|^2t^{-1-\alpha-\varepsilon}&\leq
   \limsup\limits_{t\rightarrow\infty}(\sup\limits_{n(t)-1\leq s\leq n(t)}|x(s)|^2(n(t)-1)^{-1-\alpha-\varepsilon})\limsup\limits_{t\rightarrow\infty}\frac{t^{-1-\alpha-\varepsilon}}{(n(t)-1)^{-1-\alpha-\varepsilon}})
  \\
  &=\limsup\limits_{t\rightarrow\infty}(\sup\limits_{n-1\leq s\leq n}|x(s)|^2(n-1)^{-1-\alpha-\varepsilon})\limsup\limits_{t\rightarrow\infty}\frac{t^{-1-\alpha-\varepsilon}}{(n(t)-1)^{-1-\alpha-\varepsilon}})~~~a.s.
  \\
  &=0.
   \end{split}
\end{equation*}
Due to the arbitrary of $\varepsilon$, this can imply
\begin{equation*}
  \limsup\limits_{t\rightarrow\infty}\frac{log|x(t)|}{logt}\leq \frac{1+\alpha}{2}.
\end{equation*}
\qed
\section{Consistence and convergence of the semi-implicit method }
In this section, we will employ the semi-implicit Euler methods to solve the equation (\ref{3.1}). We define a family of meshes with fixed step-size on the interval $[0,T]$, i.e.
\begin{equation}\label{4.1}
  T_N=\{t_0,t_1,t_2,\cdots,t_N\},~~~~~t_n=nh,~n=0,1,2,\cdots,N,~h=\frac{T}{N}<1.
\end{equation}
Since the points $qt_n$ will probably not be included in $T_N$, so we need another non-uniform mesh which consists of all the points $t_n$ and $qt_n$. Let
\begin{equation}\label{4.2}
  S_{N'}=\{0=t_0=s_0,s_1,s_2,\cdots,s_{N'}=T\}.
\end{equation}
For any $l$, we have $s_l=t_n$ or $s_l=qt_m$, where $t_n, t_m\in T_N$. We can also display any $s_l\in S_{N'}$ with $t_n<s_l\leq t_{n+1}$ by
\begin{equation*}
  s_l=t_n+\zeta h,~~~~~~\zeta\in (0,1], ~~t_n,t_{n+1}\in T_N.
\end{equation*}

In this paper, we denote by $y(t_n)$ the approximation of $x(t_n)$ at the point $t_n\in T_N$, and $y(qt_n)$ the approximation of $x(qt_n)$ at the point $qt_n\in S_{N'}$, then the semi-implicit Euler method is given by
\begin{equation}\label{4.3}
  y(t_{n+1})=y(t_n)+h[(1-\theta)f(t_n,y(t_n),y(qt_n))+\theta f(t_{n+1},y(t_{n+1}),y(qt_{n+1}))]+g(t_n,y(t_n),y(qt_n))\triangle B_n,
\end{equation}
where $y(t_0)=x_0$, $n=0,1,2,\cdots,N-1$, $\triangle B_n=B(t_{n+1})-B(t_n)$, $\theta\in [0,1]$. Here we require $y(t_n)$ to be  $\mathcal{F}_{t_n}$-measurable at the point $t_n,n=0,1,\cdots,N$.

We can also express (\ref{4.3}) equivalently as
\begin{equation}\label{4.4}
\begin{split}
 y(t_{n+1})=y(t_n)&+\int_{t_n}^{t_{n+1}}[(1-\theta)f(t_n,y(t_n),y(qt_n))+\theta f(t_{n+1},y(t_{n+1}),y(qt_{n+1}))]{\rm d}t
 \\
&+\int_{t_n}^{t_{n+1}}g(t_n,y(t_n),y(qt_n)){\rm d}B(t).
\end{split}
\end{equation}

Note that we can't express $y(s_l)$ which equals to $y(qt_{m})$ in(\ref{4.4}), $s_l\in S_{N'}$, so we need a continuous extension that permits the evaluation of $y(s_l)$at any point $s_l=t_n+\zeta h\in S_{N'}$, $\zeta\in (0,1]$, so we define
\begin{equation}\label{4.5}
\begin{split}
 y(s_l)=y(t_n)&+\int_{t_n}^{t_n+\zeta h}[(1-\theta)f(t_n,y(t_n),y(qt_n))+\theta f(t_{n+1},y(t_{n+1}),y(qt_{n+1}))]{\rm d}t
 \\
&+\int_{t_n}^{t_n+\zeta h}g(t_n,y(t_n),y(qt_n)){\rm d}B(t).
\end{split}
\end{equation}

For any $\theta$ given, $t_n\in T_n$, $\zeta \in (0,1]$, the local truncation error of semi-implicit Euler method for (\ref{3.1}) can be denoted by $\delta_h(t_n,\zeta)$,
\begin{equation}\label{4.6}
\begin{split}
\delta_h(t_n,\zeta)=& x(t_n+\zeta h)-\{x(t_n)+(1-\theta)\int_{t_n}^{t_n+\zeta h}f(t_n,x(t_n),x(qt_n)){\rm d}t
\\
&+\theta\int_{t_n}^{t_n+\zeta h}f(t_{n+1},x(t_{n+1}),x(qt_{n+1})){\rm d}t
\\
&+\int_{t_n}^{t_n+\zeta h}g(t_n,x(t_n),x(qt_n)){\rm d}B(t)\}.
\end{split}
\end{equation}
\begin{definition}
(i) The semi-implicit Euler method is called to be consistent with order $p_1$ in average sense, if there exist constant $C>0$ and $p_1$, which are independent of step size $h$, such that
\begin{equation}\label{4.7}
\max\limits_{0\leq n\leq N-1}\sup\limits_{\zeta\in (0,1]}|E(\delta_h(t_n,\zeta))|\leq Ch^{p_1}~~as~h\rightarrow0.
\end{equation}

(ii) The semi-implicit Euler method is called to be consistent with order $p_2$ in the sense of mean-square, if there exist constant $C$ and $p_2$, which are independent of step size $h$, such that
\begin{equation}\label{4.8}
\max\limits_{0\leq n\leq N-1}\sup\limits_{\zeta\in (0,1]}(E(|\delta_h(t_n,\zeta)|^2))^{\frac{1}{2}}\leq Ch^{p_2}~~as~ h\rightarrow0.
\end{equation}
\end{definition}

For any $\theta$ given, $t_n\in T_n$, $\zeta \in (0,1]$, the global error of semi-implicit Euler method can be denoted by $\epsilon(s_l)$
\begin{equation}\label{4.9}
 \epsilon(s_l)=\epsilon(t_n+\zeta h)=x(t_n+\zeta h)-y(t_n+\zeta h).
\end{equation}
\begin{definition}
The semi-implicit Euler method is called to be convergent with order $p$, if there exist constant $C$ and $p$, which are independent of step size $h$, such that
\begin{equation}\label{4.10}
  \max\limits_{s_l\in S_{N'}}(E(|\epsilon(s_l)|^2))^\frac{1}{2}\leq Ch^p~~as ~h\rightarrow0.
\end{equation}
\end{definition}
\begin{lemma}\cite{wc}
Assume that there exists a positive constant $K$ such that\\
(i)(Lipschitz condition) For all $t\in [0,T]$, $x_1,x_2,y_1,y_2\in R^d$,
\begin{equation*}
  |f(t,x_1,y_1)-f(t,x_2,y_2)|^2\vee|g(t,x_1,y_1)-g(t,x_2,y_2)|^2\leq K(|x_1-x_2|^2+|y_1-y_2|^2),
\end{equation*}
(ii)(Linear growth condition) For all $(t,x,y)\in [0,T]\times R^d\times R^d$,
\begin{equation*}
  |f(t,x,y)|^2\vee|g(t,x,y)|^2\leq K(1+|x|^2+|y|^2).
\end{equation*}
Then there exists a unique solution $x(t)$ to (\ref{3.1}), and $E(\sup\limits_{0\leq t\leq T}|x(t)|^2)\leq M$.
\end{lemma}
\begin{remark}~Due to Lipschiz condition and $f(t,0,0)=0,~g(t,0,0)=0$, it is not difficult to know $|f(t,x,y)|^2\leq K(|x|^2+|y|^2)$ and $|g(t,x,y)|^2\leq K(|x|^2+|y|^2).$
\end{remark}
\begin{theorem}~Under the Lipschitz condition, the semi-implicit Euler method for equation (\ref{3.1}) is consistent
(i) with  order 1.5 in average sense;
(ii) with  order 1 in mean-square sense.
\end{theorem}
{\bf Proof.}~(i) For the equation (\ref{3.1}) and the semi-implicit method (\ref{4.3}), the local truncation error takes the special form:
 \begin{equation}\label{4.11}
\begin{split}
\delta_h(t_n,\zeta)=&(1-\theta)\int_{t_n}^{t_n+\zeta h}(f(t,x(t),x(qt))-f(t_n,x(t_n),x(qt_n))){\rm d}t
\\
&+\theta\int_{t_n}^{t_n+\zeta h}(f(t,x(t),x(qt))-f(t_{n+1},x(t_{n+1}),x(qt_{n+1}))){\rm d}t
 \\
&+\int_{t_n}^{t_n+\zeta h}(g(t,x(t),x(qt))-g(t_n,x(t_n),x(qt_n))){\rm d}B(t),
 \end{split}
\end{equation}
for $n=0,1,2,\cdots,N$. we will frequently make use of H\"{o}der inequality in the next content. Note that $E(|x|)\leq (E(|x|^2))^{\frac{1}{2}}$, so taking expectation and absolute both sides of the equation above, we can estimate
\begin{equation}\label{4.12}
  \begin{split}
  |E(\delta_h(t_n,\zeta))|\leq&(1-\theta)\int_{t_n}^{t_n+\zeta h}E(|f(t,x(t),x(qt))-f(t_n,x(t_n),x(qt_n))|){\rm d}t
\\
&+\theta\int_{t_n}^{t_n+\zeta h}E(|f(t,x(t),x(qt))-f(t_{n+1},x(t_{n+1}),x(qt_{n+1}))|){\rm d}t
\\
\leq&(1-\theta)K^{\frac{1}{2}}\int_{t_n}^{t_n+\zeta h}[E(|x(t)-x(t_n)|^2)+E(|x(qt)-x(qt_n)|^2)]^{\frac{1}{2}}{\rm d}t
\\
&+\theta K^{\frac{1}{2}}\int_{t_n}^{t_n+\zeta h}[E(|x(t)-x(t_{n+1})|^2)+E(|x(qt)-x(qt_{n+1})|^2)]^{\frac{1}{2}}{\rm d}t.
  \end{split}
\end{equation}
By Lemma \ref{4.3}, $h<1$ and the integral we can obtain
\begin{equation*}
E(|x(t)-x(t_n)|^2)\leq(2K(t-t_n)+2K)E(\int_{t_n}^t(|x(t)|^2+|x(qt)|^2){\rm d}t) \leq8KMh.
\end{equation*}

In the same way, we can compute $E(|x(qt)-x(qt_n)|^2)\leq 4KMq(1+q)h$, $E(|x(t)-x(t_{n+1})|^2)\leq 8KMh$ and  $E(|x(qt)-x(qt_{n+1})|^2)\leq 4KMq(1+q)h$. Thus
\begin{equation}\label{4.13}
  |E(\delta_h(t_n,\zeta))|\leq \int_{t_n}^{t_n+\zeta h}C_1h^{\frac{1}{2}}{\rm d}t=C_1h^{\frac{3}{2}}\zeta\leq C_1h^{\frac{3}{2}}
\end{equation}
where $C_1=K(8M+4Mq(1+q))^{\frac{1}{2}}$.

This implies
$$\max\limits_{0\leq t_n\leq N-1}\sup\limits_{\zeta\in (0,1]}|E(\delta_h(t_n,\zeta))|\leq C_1h^{\frac{3}{2}}.$$
(ii) According to the definition of $\delta_h(t_n,\zeta)$, the following inequality holds.
\begin{equation*}
 \begin{split}
E(|\delta_h(t_n,\zeta)|^2)\leq&3(1-\theta)^2E(|\int_{t_n}^{t_n+\zeta h}(f(t,x(t),x(qt))-f(t_n,x(t_n),x(qt_n))){\rm d}t|^2)\\
&+3\theta^2E(|\int_{t_n}^{t_n+\zeta h}(f(t,x(t),x(qt))-f(t_{n+1},x(t_{n+1}),x(qt_{n+1}))){\rm d}t|^2)\\
&+3E(|\int_{t_n}^{t_n+\zeta h}(g(t,x(t),x(qt))-g(t_n,x(t_n),x(qt_n))){\rm d}B(t)|^2)
\end{split}
\end{equation*}
\begin{equation}\label{4.14}
  \begin{split}
\leq&3K((1-\theta)^2\zeta h+1)E(\int_{t_n}^{t_n+\zeta h}(|x(t)-x(t_n)|^2+|x(qt)-x(qt_n)|^2){\rm d}t)
\\
&+3K\theta^2\zeta hE(\int_{t_n}^{t_n+\zeta h}(|x(t)-x(t_{n+1})|^2+|x(qt)-x(qt_{n+1})|^2){\rm d}t)
\\
\leq&C_2h^2,
\end{split}
\end{equation}
where $C_2=24K^2M(\theta^2-\theta+1)(q^2+q+2)$. Let $C_3=\sqrt{C_2}$, then
 $$\max\limits_{0\leq t_n\leq N-1}\sup\limits_{\zeta\in (0,1]}(E(|\delta_h(t_n,\zeta)|^2))^{\frac{1}{2}}\leq C_3h.$$
 \qed
\begin{theorem}~Under Lipschitz condition, the semi-implicit Euler method for problem (\ref{3.1}) is convergent with order 0.5.
\end{theorem}
{\bf Proof.} For any $s_l=t_n+\zeta h\in S_{N'}$, set
\begin{equation}\label{4.15}
\begin{split}
  \nu_{h}(t_n,\zeta)=&(1-\theta)\int_{t_n}^{t_n+\zeta h}(f(t_n,x(t_n),x(qt_n))-f(t_n,y(t_n),y(qt_n))){\rm d}t
  \\
  &+\theta\int_{t_{n}}^{t_n+\zeta h}(f(t_{n+1},x(t_{n+1}),x(qt_{n+1}))-f(t_{n+1},y(t_{n+1}),y(qt_{n+1}))){\rm d}t
  \\
  &+\int_{t_n}^{t_n+\zeta h}(g(t_n,x(t_n),x(qt_n))-g(t_n,y(t_n),y(qt_n))){\rm d}B(t).
  \\
  \end{split}
\end{equation}
Then
\begin{equation}\label{4.16}
  \epsilon(s_l)=x(t_n+\zeta h)-y(s_l)=\epsilon(t_n)+\delta_h(t_n,\zeta)+  \nu_{h}(t_n,\zeta).
\end{equation}

Squaring both sides of the equation above, employing the conditional expectation with respect to the $\sigma$-algebra $\mathcal{F}_0$, and taking absolute values, we get
\begin{equation}\label{4.17}
\begin{split}
E(|\epsilon(s_l)|^2|\mathcal{F}_0)\leq&E(|\epsilon(t_n)|^2|\mathcal{F}_0)+E(|\delta_h(t_n,\zeta)|^2|\mathcal{F}_0)+E(|\nu_{h}(t_n,\zeta)|^2|\mathcal{F}_0)
\\
&+2E(|\epsilon(t_n)|\cdot|\delta_h(t_n,\zeta)||\mathcal{F}_0)+2E(|\epsilon(t_n)|\cdot|\nu_h(t_n,\zeta)||\mathcal{F}_0)
\\
&+2E(|\delta_h(t_n,\zeta)|\cdot|\nu_h(t_n,\zeta)||\mathcal{F}_0) \quad \quad\quad\quad \quad\quad\quad \quad\quad\quad \quad \quad a.s.
\\
=&A_1+A_2+A_3+A_4+A_5+A_6.
 \end{split}
\end{equation}
Next we will estimate the six terms in (\ref{4.14}). For the term $A_2$, by (\ref{4.14}) we have
\begin{equation*}
  A_2=E(|\delta_h(t_n,\zeta)|^2|\mathcal{F}_0)=E(E(|\delta_h(t_n,\zeta)|^2|\mathcal{F}_n)|\mathcal{F}_0)\leq C_2h^2.
\end{equation*}
For $A_3$ in (\ref{4.17}), we obtain
\begin{equation*}
  \begin{split}
  A_3\leq&3(1-\theta)^2K\zeta hE(\int_{t_n}^{t_n+\zeta h}(|x(t_n)-y(t_n)|^2+|x(qt_n)-y(qt_n)|^2){\rm d}t|\mathcal{F}_0)
  \\
   &+3\theta^2K\zeta hE(\int_{t_{n}}^{t_n+\zeta h}(|x(t_{n+1})-y(t_{n+1})|^2+|x(qt_{n+1})-y(qt_{n+1})|^2{\rm d}t|\mathcal{F}_0)
   \\
   &+3KE(\int_{t_n}^{t_n+\zeta h}(|x(t_n)-y(t_n)|^2+|x(qt_n)-y(qt_n)|^2){\rm d}t|\mathcal{F}_0)
   \\
   =&3K\zeta h((1-\theta)^2\zeta h+1)E(|\epsilon(t_n)|^2|\mathcal{F}_0)+3K\zeta h((1-\theta)^2\zeta h+1)E(|\epsilon(qt_n)|^2|\mathcal{F}_0)
   \\
   &+3\theta^2K\zeta^2h^2E(|\epsilon(t_{n+1})|^2|\mathcal{F}_0)+3\theta^2K\zeta^2h^2E(|\epsilon(qt_{n+1})|^2|\mathcal{F}_0).
  \end{split}
\end{equation*}
We estimate $A_4$,
\begin{equation*}
\begin{split}
A_4\leq 2(E(|E(\delta_h(t_n,\zeta))|^2|\mathcal{F}_0))^\frac{1}{2}\cdot (E(|\epsilon(t_n)|^2|\mathcal{F}_0))^\frac{1}{2}
\leq C_1^{2}h^{2}+hE(|\epsilon(t_n)|^2|\mathcal{F}_0).
\end{split}
\end{equation*}
In the same way, we can see
\begin{equation*}
\begin{split}
  A_5\leq &(2-\theta)K^{\frac{1}{2}}\zeta hE(|\epsilon(t_n)|^2|\mathcal{F}_0)+(1-\theta)K^{\frac{1}{2}}\zeta hE(|\epsilon(qt_n)|^2|\mathcal{F}_0)
  \\
  &+\theta K^{\frac{1}{2}}\zeta hE(|\epsilon(t_{n+1})|^2|\mathcal{F}_0)+\theta K^{\frac{1}{2}}\zeta hE(|\epsilon(qt_{n+1})|^2|\mathcal{F}_0),
  \end{split}
\end{equation*}
and
\begin{equation*}
    \begin{split}
 A_6\leq&2(E(|\nu_h(t_n,\zeta)|^2|\mathcal{F}_0))^{\frac{1}{2}}\cdot(E(|\delta_h(t_n,\zeta)|^2|\mathcal{F}_0))^{\frac{1}{2}}
    \\
    \leq&C_2h^2+3K\zeta h((1-\theta)^2\zeta h+1)E(|\epsilon(t_n)|^2|\mathcal{F}_0)
    \\
    &+3\theta^2K\zeta^2 h^2E(|\epsilon(t_{n+1})|^2|\mathcal{F}_0)
    \\
    &+3K\zeta h((1-\theta)^2\zeta h+1)E(|\epsilon(qt_n)|^2|\mathcal{F}_0)
    \\
    &+3\theta^2K\zeta^2 h^2E(|\epsilon(qt_{n+1})|^2|\mathcal{F}_0).
    \end{split}
\end{equation*}

Combining these results, we can compute
\begin{equation}\label{4.18}
  \begin{split}
  E(|\epsilon(t_n+\zeta h)|^2|\mathcal{F}_0)\leq& (1+6K\zeta h((1-\theta)^2\zeta h+1)+h+(2-\theta)K^{\frac{1}{2}}\zeta h)E(|\epsilon(t_n)|^2|\mathcal{F}_0)
  \\
  &+(6\theta^2K\zeta^2h^2+\theta K^{\frac{1}{2}}\zeta h)E(|\epsilon(t_{n+1})|^2|\mathcal{F}_0)
  \\
  &+(6K\zeta h((1-\theta)^2\zeta h+1)+(1-\theta)K^{\frac{1}{2}}\zeta h)E(|\epsilon(qt_n)|^2|\mathcal{F}_0)
  \\
  &+(6\theta^2K\zeta^2h^2+\theta K^{\frac{1}{2}}\zeta h)E(|\epsilon(qt_{n+1})|^2|\mathcal{F}_0)+(2C_2+C_1^2)h^2.
  \end{split}
\end{equation}

Set $R_0=0$, $R_n=\max\limits_{0\leq i< n}\sup\limits_{\zeta\in (0,1]}E(|\epsilon(t_i+\zeta h)|^2|\mathcal{F}_0)$, then
$$E(|\epsilon(t_n)|^2|\mathcal{F}_0)\leq R_n, ~~ E(|\epsilon(qt_n)|^2|\mathcal{F}_0)\leq R_n.$$
In (\ref{4.18}), we need to calculate $E(|\epsilon(qt_{n+1})|^2|\mathcal{F}_0)$, which depends on either $t_n<qt_{n+1}<t_{n+1}$ or $qt_{n+1}<t_n$.

Case 1: If $t_n<qt_{n+1}<t_{n+1}$, then $E(|\epsilon(qt_{n+1})|^2|\mathcal{F}_0)\leq R_{n+1}$. According to (\ref{4.18}), we can see
\begin{equation}
  \begin{split}
E(|\epsilon(t_n+\zeta h)|^2|\mathcal{F}_0)\leq&(1+12(1-\theta)^2Kh^2+12Kh+h+(3-2\theta)K^{\frac{1}{2}}h)R_n\\
  &+(12\theta^2Kh^2+2\theta K^{\frac{1}{2}}h)R_{n+1}.
\end{split}
\end{equation}
So
\begin{equation*}
\begin{split}
  R_{n+1}=&\max\limits_{0\leq i< n+1}\sup\limits_{\zeta\in (0,1]}E(|\epsilon(t_i+\zeta h)|^2|\mathcal{F}_0)
  \\
  \leq&(1+12(1-\theta)^2Kh^2+12Kh+h+(3-2\theta)K^{\frac{1}{2}}h)R_n
  \\
  &+(12\theta^2Kh^2+2\theta K^{\frac{1}{2}}h)R_{n+1}+(2C_2+C_1^2)h^2.
  \end{split}
\end{equation*}
There is $h_0=\frac{\sqrt{13}-1}{12}K^{-\frac{1}{2}}$, such that $1-12\theta^2Kh^2-2\theta K^{\frac{1}{2}}h>0$ when $0<h<h_0$. Therefore
\begin{equation*}
  R_{n+1}\leq(1+h\frac{1+12(1-\theta)^2K+12K+3K^{\frac{1}{2}}+12\theta^2K}{1-12\theta^2Kh^2-2\theta K^{\frac{1}{2}}h})R_n+\frac{2C_2+C_1^2}{1-12\theta^2Kh^2-2\theta K^{\frac{1}{2}}h}h^2.
\end{equation*}
Case 2: If $qt_{n+1}<t_n$, then $E(|\epsilon(qt_{n+1})|^2|\mathcal{F}_0)\leq R_n$. In the same way as case 1, we can get
\begin{equation*}
  R_{n+1}\leq(1+h\frac{1+12(1-\theta)^2K+12K+3K^{\frac{1}{2}}+12\theta^2K}{1-6\theta^2Kh^2-\theta K^{\frac{1}{2}}h})R_n+\frac{2C_2+C_1^2}{1-6\theta^2Kh^2-\theta K^{\frac{1}{2}}h}h^2,
\end{equation*}
when $0<h<h_1=\frac{1}{3}K^{-\frac{1}{2}}.$ Now take $0<L<1$, which is independent of $h$, such that $12\theta^2Kh^2+2\theta K^{\frac{1}{2}}h<L$, and set
$$M(\theta)=\frac{1+12(1-\theta)^2K+12K+3K^{\frac{1}{2}}+12\theta^2K}{1-L} ,~~C(\theta)=\frac{2C_2+C_1^2}{1-L}.$$

Then combining case 1 and case 2, the $R_{n+1}$ satisfies
\begin{equation*}
  \begin{split}
  R_{n+1}\leq&(1+hM(\theta))R_n+C(\theta)h^2\leq(1+hM(\theta))R_{n-1}+(1+hM(\theta))C(\theta)h^2+C(\theta)h^2
  \\
  \leq&\cdots\leq(1+hM(\theta))^{n+1}R_0+C(\theta)h^2\sum_{i=0}^{n}(1+hM(\theta))^{i}\leq \frac{(1+hM(\theta))^{n+1}-1}{M(\theta)}C(\theta)h.
  \end{split}
\end{equation*}
The expression above indicates that
\begin{equation*}
  E(|\epsilon(s_l)|^2|\mathcal{F}_0)\leq R_{n+1}\leq \frac{e^{TM(\theta)}-1}{M(\theta)}C(\theta)h,
\end{equation*}
for any $s_l=t_n+\zeta h\in S_{N'}$ holds, where $t_n\in T_N$, $\zeta\in (0,1]$.
By the definition of convergence, we can show
\begin{equation*}
  \max\limits_{s_l\in S_{N'}}(E(|\epsilon(s_l)|^2|\mathcal{F}_0))^{\frac{1}{2}}\leq \sqrt{\frac{C(\theta)(e^{TM(\theta)}-1)}{M(\theta)}}h^{\frac{1}{2}}
\end{equation*}
The theorem is proved.

\qed\\


\begin{thebibliography}{99}
\bibitem{jg1}
Appleby, J.A.D., Berkolaiko, G., Rodkina, A., Non-exponential stability and decay rates in nonlinear
stochastic difference equations with unbounded noise, Stochastics: An International Journal of Probability and Stochastics Processes, 81(2009), 99-127.
\bibitem{eb}
Appleby, J.A.D., Buckwar, E., Sufficient conditions for polynomial asymptotic behaviour of the stochastic pantograph equation, Stochastic Anal, 2003.
\bibitem{jg2}
Appleby, J.A.D., Mackey, D., Almost sure polynomial asymptotic stability of stochastic difference equations, Journal of Mathematical Sciences, 149(2008), 1629-1647.
\bibitem{jg3}
Appleby, J.A.D., Mackey, D., Polynomial Asymptotic Stability of Damped Stochastic Differential Equations, Electronic Journal Qualitative Therory of Differntial Equations. 2(2004), 1-33.
\bibitem{bb}
Baker, C.T.H., Buckwar, E., Continuous $\theta$-Methods for the Stochastic Pantograph Equation, Electronic Transactions on Numerical Analysis, 11(2000), 131-151.
\bibitem{cd}
Carr, J., Dyson, J., The functional differential equation $y'(x)=ay(\lambda x)+by(x)$, Proc. Roy. Soc. Edinburgh Sect. A. 74(1974), 165-174.
\bibitem{wc}
Fan, Z.C., Liu, M.Z., Cao, W.R., Existence and uniqueness of the solutions and convergence of semi-implicit Euler methods for stochastic pantograph equations, J. Math. Anal. Appl. 325 (2007) 1142-1159.
\bibitem{fan}
Fan, Z.C., Song, M.H., Liu, M.Z., The $\alpha$th moment stability for the stochastic pantograph equation,Journal of Computational and Applied Mathematics, 233(2009), 109-120.
\bibitem{kx}
Liu, K., Mao, X.R., Large time decay behavior of dynamical equations with random perturbation features, Stochastic Analysis and Applications, 19(2001), 295-327.
\bibitem{mz}
Liu, M.Z., Yang, Z.W., G.D.Hu, Asymptotical stability of the numerical methods with the constant stepsize for the pantograph equation, BIT, 45(2005), 743-759.
\bibitem{mw}
Mao, W., Convergence analysis of semi-implicit Euler methods for solving stochastic equations with variable delays and random Jump magnitudes, Journal of Computational and Applied Mathematics, 235(2011), 2569-2580.

\bibitem{ab}
Tsoi, A.H., Zhang, B., Weak exponential stability of stochastic differential equations, Stochastic Analysis and Applications, 15(1997), 643-649.

\end{thebibliography}
\end{document}